\documentclass[paper=a4, 12pt]{article} % A4 paper and 12pt font size

\usepackage[latin1]{inputenc} 
\usepackage{geometry}
\usepackage{amsmath,amsfonts,amsthm,amssymb,mathtools} % Math packages
\usepackage{enumerate}
\usepackage{stmaryrd}
\usepackage[english]{babel}
\usepackage{relsize}
\usepackage{color}
\usepackage{fancyhdr}
 \usepackage{multicol}
 \usepackage{tikz-cd}
 \usepackage[noadjust]{cite}
 \usepackage[hyperfootnotes=false]{hyperref}
\usepackage{url}
\numberwithin{figure}{section} % Number figures within sections (i.e. 1.1, 1.2, 2.1, 2.2 instead of 1, 2, 3, 4)
\numberwithin{table}{section} % Number tables within sections (i.e. 1.1, 1.2, 2.1, 2.2 instead of 1, 2, 3, 4)

\theoremstyle{theorem}
\newtheorem{theorem}{Theorem}[section]
\newtheorem{proposition}[theorem]{Proposition}

\newtheorem{lemma}[theorem]{Lemma}

\theoremstyle{definition}
\newtheorem{definition}[theorem]{Definition}

\newcommand{\dcorner}{\rotatebox[origin=c]{-45}{$\lrcorner$}}
\newcommand{\dpbk}{\arrow[dd,phantom,"\dcorner", very near start]}
\newcommand{\rdpbk}{\arrow[rd,phantom,"\lrcorner ", very near start]}
\newcommand{\lupbk}{\arrow[lu,phantom,"\ulcorner ", very near start]}

\DeclareMathAlphabet{\mathbbe}{U}{bbold}{m}{n}
\newcommand{\simplexcategory}{\mathbbe{\Delta}} 

\newcommand{\pow}[1]{\llbracket #1\rrbracket}

\newcommand{\infgrpd}{\mathbf{Grpd}}
\newcommand{\set}{\mathbf{Set}}

\newcommand{\sur}{\mathbf{S}}

\newcommand{\grpdSprod}{\odot}

\newcommand{\arxiv}[1]{\href{http://arxiv.org/pdf/#1}{arXiv:#1}}

\newcommand{\Q}{\mathbb Q}
\newcommand{\name}[1]{\ulcorner #1\urcorner}

\newcommand{\C}{\mathcal C}
\newcommand{\B}{\mathbb B}
\newcommand{\Sg}{\mathbb S}

\newcommand{\bbP}{\mathbb P}

\newcommand{\x}{\mathbf x}
\renewcommand{\l}{\lambda}

\DeclareMathOperator{\autiv}{aut}
\DeclareMathOperator{\aut}{Aut}

\DeclareMathOperator{\iso}{Iso}

\DeclareMathOperator{\im}{Im}
\DeclareMathOperator{\fun}{Fun}

\begin{document}
 %%----------------------------------------------------------------------------------------
%%	TITLE
%%----------------------------------------------------------------------------------------
 
 \title{Combinatorics and simplicial groupoids}
 \date{}
 \author{Alex Cebrian \\Universitat Aut\`onoma de Barcelona}
 
 \maketitle
 
 %%----------------------------------------------------------------------------------------
%%	ABSTRACT
%%----------------------------------------------------------------------------------------
 
 \begin{abstract}
This expository paper starts with a brief  survey on the relation between partitions and surjections of sets, and then gives a quick introduction to the theories of incidence algebras, Segal groupoids and combinatorial species. The aim is to explain an objective construction, in terms of simplicial groupoids, of both the Fa\`a di Bruno bialgebra and the plethystic bialgebra. This paper can be seen as an extended introduction to the author's paper \cite{Cebrian}.
 \end{abstract}
 
 \vspace{-30pt}
 %\tableofcontents

 %%----------------------------------------------------------------------------------------
%%	FOOTNOTE
%%----------------------------------------------------------------------------------------
 {\let\thefootnote\relax\footnote{
 %\textit{Date}: November 25, 2019. \\
 The  author was supported by the project FEDER-MEC MTM2016-80439-P.}}
 
 %%----------------------------------------------------------------------------------------
%%	SECTION: partitions and surjection
%%----------------------------------------------------------------------------------------
 \section{Introduction}
  It is well appreciated in combinatorics that working with the combinatorial structures themselves gives a deeper understanding than working with their numbers. This is called objective combinatorics. The theory of species, developed by Joyal \cite{Joyal:1981}, is a cornerstone in this context. It makes use of category theory to objectify generating functions of combinatorial structures.
On the other hand, decomposition spaces (certain simplicial spaces) provide a general machinery to objectify the notion of incidence algebra in algebraic combinatorics. Decomposition spaces were introduced by  G\'alvez, Kock and Tonks \cite{GKT:HLA,GKT:DSIAMI-1,GKT:DSIAMI-2} in this framework and they are the same as $2$-Segal spaces, introduced by Dyckerhoff and Kapranov \cite{DK} in the context of homological algebra and representation theory.

 In Section \ref{preliminaries} we introduce incidence coalgebras and explain how can they be encoded in Segal groupoids, which are a particular simpler case of decomposition space. We apply this to express two incidence coalgebras in terms of Segal groupoids: the Fa\`a di Bruno bialgebra and the plethystic bialgebra. The first one comes from substitution of power series in one variable, and it is closely related to the theory of species. We explain this in Section \ref{species}. The second one comes from plethystic substitution of power series in infinitely many variables, and it is closely related to the theory of partitionals \cite{Nava-Rota}, a generalization of species to partitions. We explain this in Section \ref{partitionals}. Both the theory of species and the theory of partitionals are based on sets and partitions. However, to give an interpretation of these bialgebras through simplicial groupoids it is better to work with the equivalent category of sets and surjections, as advocated in \cite{GKT:FdB} and \cite{Cebrian}. Section \ref{partitions} is a survey on the relation between partitions and surjections. The two Segal goupoids realizing  the Fa\`a di Bruno bialgebra and the plethystic bialgebra are respectively $N\sur$ (Section \ref{species}), the fat nerve of the category of sets and surjections, and $T\sur$ (Section \ref{partitionals}), introduced in \cite{Cebrian}.

%%----------------------------------------------------------------------------------------
%%	SECTION: partitions and surjection
%%----------------------------------------------------------------------------------------

\section{Partitions and surjections}\label{partitions}

The content on partitions featuring in this section has been taken from \cite{Nava-Rota}. A \emph{partition} $\pi$ of a finite set $E$ is a family of subsets of $E$, called \emph{blocks}, such that every block of $\pi$ is nonempty, the blocks are pairwise disjoint, and every element of $E$ is contained in some block. Given two partitions $\pi,\sigma$ of $E$ we say that $\pi$ is finer than $\sigma$ (or $\sigma$ is coarser than $\pi$), and denote by $\pi\leq\sigma$ if every block of $\pi$ is a subset of some block of $\sigma$. In this case we define the \emph{induced partition} $\sigma | \pi$ to be the partition on the set of blocks of $\pi$ given by the blocks of $\sigma$. Also, the \emph{restriction} of a partition $\pi$ to a subset $B\subseteq E$ is the partition of $B$ given by the intersections of $B$ with the blocks of $\pi$ and denoted by $\pi_B$.

Notice that the relation $\leq$ defines a partial order on the set $\Pi(E)$ of all partitions of $E$, and this order has a minimum, given by the partition with unitary blocks and denoted by $\hat{0}$, and a maximum, given by the partition with one block and denoted by $\hat{1}$. 
In particular, for every $\pi,\sigma\in \Pi(E)$  the supremum $\pi\vee\sigma$ and the infimum $\pi\wedge\sigma$ of $\pi$ and $\sigma$ exist in $(\Pi(E),\leq)$, and they are respectively called the \emph{join} and the \emph{meet}. For example if $E=\{1,2,3,4,5,6\}$, $\pi=\{\{1,2\},\{3,4,5\},\{6\}\}$ and $\sigma=\{\{1,2,6\},\{3,4\},\{5\}\}$ then
$$\pi\wedge\sigma=\{\{1,2\},\{3,4\},\{5\},\{6\}\}, \;\;\;\pi\vee\sigma=\{\{1,2,6\},\{3,4,5\}\},$$
$$ \sigma|(\pi\vee\sigma)=\{\{\{1,2,6\}\},\{\{3,4\},\{5\}\}\} \;\;\;\text{and}\;\;\; \pi_{\{1,3,4,6\}}=\{\{1\},\{3,4\},\{6\}\}.$$
It is easy to see that in fact
$$\pi\wedge\sigma=\{B\cap C \;|\; B\in \pi, C\in \sigma, B\cap C\neq \o\}.$$
The join is, roughly speaking, the union of all blocks with common elements. However to give a precise and simple definition of the join it is preferable to view partitions as equivalence relations. This will also help us to introduce other notions later. It is clear that a partition $\pi$ on a set $E$ defines an equivalence relation $\sim_{\pi}$ on $E$, where two elements are related if they belong to the same block of $\pi$. In this setting the meet of two partitions $\pi,\sigma$ is given by (for $p,q\in E$) $p\sim_{\pi\wedge\sigma}q$ if $p\sim_{\pi}q$ and $p\sim_{\sigma}q$, and their join is given by $p\sim_{\pi\vee\sigma}q$ if there exists a finite sequence $r_0,\dots,r_n$ such that 
\begin{equation}\label{joinrelation}
p=r_0\sim_1 r_1\sim_2\dots \sim_n r_n=q,
\end{equation}
 where each relation $\sim_k$ is either $\sim_{\pi}$ or $\sim_{\sigma}$. We say that two partitions $\pi,\sigma$ are \emph{independent} if every block of $\pi$ meets every block of $\sigma$. We say that they \emph{commute} if, for every $p,q\in E$, we have that $p\sim_{\pi}r \sim_{\sigma}q$ for some $r\in E$ if and only if $p\sim_{\sigma}s \sim_{\pi}q$ for some $s\in E$. It is straightforward to see that two partitions are independent if and only if they commute and their join is $\hat{1}$. The following result says that commuting is the same as being blockwise independent.
 \begin{proposition} \label{comind}Let $\pi,\sigma\in \Pi(E)$. Then $\pi$ and $\sigma$ commute if and only if or every $B\in \pi\vee\sigma$ the restrictions $\pi_B$ and $\sigma_B$ are independent partitions of $B$.
 \end{proposition}
 We are now ready define the most intricate and important definition regarding partitions in this survey. \begin{definition}[\cite{Nava-Rota}] Let $\sigma$ be a partition of $E$. A pair $(\pi,\tau)$ of partitions of $E$ is called a \emph{transversal} of $\sigma$ when
 \begin{multicols}{2}
 \begin{enumerate}[(i)]
 \item $\pi\leq \sigma$,
 \item $\pi\wedge \tau=\hat{0}$,
 \item $\pi$ and $\tau$ commute, and
 \item $\sigma\vee\tau=\pi\vee\tau$.
 \end{enumerate}
 \end{multicols}
 \end{definition}

There is a category $\bbP$ whose objects are pairs $(E,\pi)$, where $E$ is a set and $\pi\in \Pi(E)$, and whose morphisms are defined as follows: if $(F,\sigma)$ is another object of $\bbP$ a morphism
$$f\colon (E,\pi)\longrightarrow(F,\sigma)$$
is a bijection $f\colon E\rightarrow F$ which maps blocks of $\pi$ to blocks of $\sigma$. Notice that $\bbP$ is in fact a groupoid, since all the morphisms are invertible. In this groupoid, the isomorphism class of a partition $(E,\pi)$ can be described by the sequence of natural numbers
$$\l=(\l_1,\l_2,\dots), \;\; \text{where } \l_k=\text{number of blocks of size }k\text{ of }(E,\pi).$$
Observe that $|E|=1\cdot \l_1+2\cdot \l_2+3\cdot \l_3+\cdots$ and that the number of blocks of $\pi$ is $|\pi|=\l_1+\l_2+\cdots$. Also, notice that the number of automorphisms of $(E,\pi)$ is
$$\autiv(\l)=1!^{\l_1}\l_1!\cdot2!^{\l_2}\l_2!\cdot 3!^{\l_3}\l_3!\cdots,$$
because an automorphism of $\pi$ permutes the elements inside each block and permutes the blocks of the same size.

Let us translate all these notions from the category of partitions $\bbP$ to the equivalent category of surjection $\Sg$. Among the advantages of surjections over partitions there is the fact that partitions of partitions are pairs of composable surjections. The category of surjections has as objects surjections $E\twoheadrightarrow S$ between finite sets and as morphisms commutative squares
\begin{equation*}
\begin{tikzcd}[sep={3em,between origins}]
E\arrow[d,twoheadrightarrow]\arrow[r,"\sim"]&F\arrow[d,twoheadrightarrow]\\
S\arrow[r,"\sim"]&R,
\end{tikzcd}
\end{equation*}
where the horizontal arrows are bijections. It is clear that $f\colon E\twoheadrightarrow S$ corresponds to the partition $\pi$ of $E$ given by $p\sim_{\pi}q$ if $f(p)=f(q)$ or, what is the same, the partition whose blocks are the fibres of  $f$. Hence, the isomorphism class of a surjection is given by a sequence $\l=(\l_1,\l_2,\dots)$ where $\l_k$ is the number of fibres of size $k$, and the number of automorphisms of $f$ is also $\autiv(\l)$: in this case $\l_1!\cdot\l_2!\cdots$ is the number of bijections $S\xrightarrow{\sim}S$ permuting elements with a fibre of the same size, and $1!^{\l_1}\cdot 2!^{\l_2}\cdots$ is the number of fibrewise bijections $E\xrightarrow{\sim}E$.

Consider $\pi,\tau\in \Pi(E)$ and let $\pi\colon E\twoheadrightarrow S$  and  $\tau\colon E\twoheadrightarrow X$ be their corresponding surjections. Construct the diagram of sets
\begin{equation*}\label{partsurj}
\begin{tikzcd}[column sep=small]
E\arrow[dr,dashrightarrow,"\phi" description]\arrow[drr,twoheadrightarrow,bend left,"\tau"]\arrow[ddr,twoheadrightarrow,bend right,"\pi"']&&\\
&S{\times_I}X\arrow[d,twoheadrightarrow]\arrow[r,twoheadrightarrow]\rdpbk&X\arrow[d,twoheadrightarrow]\\
&S\arrow[r,twoheadrightarrow]&I \lupbk
\end{tikzcd}
\end{equation*}
by taking pushout along $\pi$ and $\tau$ and pullback of the pushout diagram. Note that all the arrows are surjections except perhaps $\phi$. Note also that any pullback of surjections is also a pushout square.  \begin{lemma} \label{lema:PS}Let $\pi$ and $\tau$ be two partitions of $E$ presented as surjections as above, and  $\sigma\colon E\twoheadrightarrow B$ another partition. Let also $A\subseteq E$. 
  \begin{enumerate}[(i)]
  \item $\pi\leq \sigma$ if and only if $\sigma$ factors through $\pi$: $E\twoheadrightarrow S\twoheadrightarrow B$. Moreover the surjection $S\twoheadrightarrow B$ corresponds to $\sigma|\pi$.
  \item $\pi_A$ corresponds to the unique surjection $A\twoheadrightarrow R$ (up to isomorphism) that factors the morphism $A\hookrightarrow E\twoheadrightarrow S$ as a surjection followed by an injection $A\twoheadrightarrow R\hookrightarrow S$.
  \item $\hat{0}$ is $E\twoheadrightarrow E$ and $\hat{1}$ is $E\twoheadrightarrow 1$.
  \item The join $\pi\vee \tau$ corresponds to the pushout surjection $E\twoheadrightarrow I$. 
  \item The meet $\pi \wedge \tau$ corresponds to the surjection $\phi\colon E\twoheadrightarrow \im (\phi)$. Hence $\pi \wedge \tau=\hat{0}$ if and only if $\phi$ is injective. \item $\pi$ and $\tau$ commute if and only if $\phi$ is surjective.
 \item $\pi$ and $\tau$ are independent if and only if $\phi$ is surjective and $I=1$.
 \end{enumerate}
 \begin{proof} (i), (ii) and (iii) are clear. (iv) follows from the fact that the pushout is precisely $X\sqcup_{/\sim} S$ where $\sim$ is the equivalence relation generated by the relation of belonging to the same fibre along $\pi$ and $\tau$. This is precisely the same relation defined in (\ref{joinrelation}). 
 
 For (v), recall that for every $p,q\in E$ we have that $p\sim_{\pi\wedge\tau}q$ if and only if $p\sim_{\pi}q$ and $p\sim_{\tau}q$, but this is the same as $\pi(p)=\pi(q)$ and $\tau(p)=\tau(q)$, which is the same as $\phi(p)=\phi(q)$. But this is equivalent to $p\sim_{\phi}q$, considering $\phi$ as a surjection to its image. 
 
Finally, if $\pi$ and $\tau$ are independent then $\pi\wedge\tau=\hat{1}$, so that $I=1$, and every fibre along $\pi$ has nonempty intersection with every fibre along  $\tau$, which means that $\phi$ is surjective. The converse is similar. This, together with Proposition \ref{comind} shows (vi), since the set $\pi\vee\tau$ is precisely $I$.
 \end{proof}
 \end{lemma}
 In view of this lemma, a transversal of the surjection  $\sigma\colon E\twoheadrightarrow B$ is a diagram 
\begin{equation}
\begin{tikzcd}[sep={3em,between origins}]
& & E\arrow[dl,twoheadrightarrow,"\pi"']\arrow[dr,twoheadrightarrow,"\tau"]\arrow[ddll,twoheadrightarrow, bend right,"\sigma"']\dpbk& &\\
&S\arrow[dl,twoheadrightarrow]\arrow[dr,twoheadrightarrow]& &X\arrow[dl,twoheadrightarrow]\arrow[dr,twoheadrightarrow]&\\
B\arrow[rr,twoheadrightarrow]&&I\arrow[rr,twoheadrightarrow]&&1
\end{tikzcd}
\end{equation}
where the square is obtained as the pushout of $\pi$ and $\tau$. The fact that $\pi\wedge \tau=\hat{0}$ and that $\pi$ and $\tau$ commute implies that this square is also a pullback. Furthermore the condition that the pushouts $\pi\vee\tau$ and $\sigma\vee \tau$ coincide gives a map $B\twoheadrightarrow I$. Conversely, any commutative square of the form
\begin{equation*}
\begin{tikzcd}[sep={3em,between origins}]
&S\arrow[dl,twoheadrightarrow]\arrow[dr,twoheadrightarrow]&\\
B\arrow[dr,twoheadrightarrow]&&I\\
&I\arrow[ur,equal]&
\end{tikzcd}
\end{equation*}
is a pushout in the category of surjections. Therefore the  map $B\twoheadrightarrow I$ says that $\pi\vee\tau$  coincides with $\sigma\vee \tau$. This diagramatic rendition of the notion of transversal \cite{Cebrian} will be the key to Section \ref{partitionals}.

 %%----------------------------------------------------------------------------------------
%%	SECTION: incidence coalgebras Segal groupoids
%%----------------------------------------------------------------------------------------
 
\section{Incidence coalgebras and Segal groupoids}\label{preliminaries}\label{HoCard}

Coalgebras arise in algebraic combinatorics from the ability to decompose structures into smaller ones. Recall that a \emph{coalgebra} is the dual notion of a unital associative algebra. That is, a vector space $V$ over a field $k$ together with $k$-linear maps $\Delta\colon V\rightarrow V\otimes V$ and $\epsilon\colon V\rightarrow k$, called comultiplication and counit respectively, satisfying 
$$(id_V\otimes \Delta)\circ \Delta=(\Delta\otimes id_V)\circ \Delta\;\;\; \text{ and }\;\;\; (id_V\otimes \epsilon)\circ \Delta=id_V=(\epsilon\otimes id_V)\circ \Delta.$$
In combinatorics, these vector spaces are usually vector spaces spanned by isomorphism classes of structures.

Rota \cite{Rota} showed  that many of these coalgebras admit an interpretation in terms of \emph{incidence coalgebras of posets}: from any locally finite poset, form the free vector space on its intervals, and endow this with a coalgebra structure by defining the comultiplication as
$$\Delta([x,y])=\sum_{x\leq m \leq y} [x,m]\otimes [m,y].$$
 Observe that a poset can be viewed as a category where there is at most one arrow between any two objects. The theory for locally finite posets was generalized to categories by Leroux \cite{Leroux}, and goes as follows: a category is \emph{locally finite} if every arrow admits only a finite number of $2$-step factorizations. The \emph{incidence coalgebra} of a locally finite category is the free vector space spanned by its arrows, with comultiplication
 $$\Delta(f)=\sum_{b\circ a=f} a\otimes b$$
 and counit  $\epsilon(id_x)=1$ and $\epsilon(f)=0$ else. The coassociativity of $\Delta$ comes from the associativity of composition of arrows. 
 
 It is well appreciated in combinatorics that bijective proofs represent a deeper insight into combinatorial structures than algebraic proofs. Lawvere and Menni \cite{LawvereMenni} pioneered the so-called objective method in this context, with the aim to work directly with the combinatorial structures, rather than their numbers, by using linear algebra over sets. Let $\set$ be the category of sets. Then the objective counterpart of the vector space spanned by a set $S$ is the slice category $\set_{/S}$ (cf. \cite{GKT:HLA}). An object in this category is a morphism $A\xrightarrow{f} S$ of sets, and it corresponds to the vector whose $s$-entry (for $s\in S$) is $|f^{-1}(s)|$. Linear maps $\set_{/S}\rightarrow \set_{/R}$ at the objective level are given by spans $S\leftarrow M \rightarrow R$, and obtained by taking pullback and postcomposition, as in (\ref{segalcom}). A coalgebra in $\set_{/S}$ is thus given by a comultiplication span $S\leftarrow M\rightarrow S\times S$  and a counit span $S\leftarrow N \rightarrow 1$. However, combinatorial structures have symmetries, and to deal with them it is useful to update this objective method to groupoids and homotopy linear algebra over groupoids. A brief introduction to the homotopy approach to groupoids in combinatorics can be found in \cite[\S 3]{GKT:FdB}. We explain the basic notions next.
  
 A \emph{groupoid} is a category whose arrows are all isomorphisms.  
 We will denote by $\infgrpd$ the category of groupoids, whose objects are groupoids and whose morphisms are functors.  An \emph{equivalence} of groupoids is an equivalence of categories. Given a groupoid $X$, we denote by $\pi_0(X)$ the set of isomorphism classes of $X$, and for $x\in X$ we denote by  $\aut(x)$ the group of automorphisms of $x$. A concrete example of this can be found in Section \ref{partitions} with the groupoids of partitions and surjections.  A groupoid $X$ is \emph{finite} if $\pi_0(X)$ is a finite set and $\aut(x)$ is a finite group for every element $x$. If only the latter is satisfied then it is called \emph{locally finite}. A morphism of groupoids is called \emph{finite} when all its homotopy fibres are finite.

 The \emph{homotopy pullback} of a diagram of groupoids $X\xrightarrow{f} B\xleftarrow{g} Y$ is the groupoid
 $Z$ whose objects are triples $(x,y,\phi)$ with $x\in X$, $y\in Y$ and $\phi\colon f(x)\rightarrow g(y)$ an arrow of $B$, and whose arrows are pairs $(\alpha,\beta)\colon (x,y,\phi)\rightarrow (x',y',\phi')$ consisting of two arrows $\alpha\colon x\rightarrow x'$ and $\beta\colon y\rightarrow y'$ satisfying $g(\beta)\circ \phi=\phi'\circ f(\alpha)\colon f(x)\rightarrow g(y')$. Given a morphism of groupoids $X\xrightarrow{f} B$ and an object $b\in B$, the \emph{homotopy fibre} of $b$ along $f$ is the groupoid $X_b$ obtained by taking the homotopy pullback of the diagram $X\xrightarrow{f} B\xleftarrow{\name{b}} 1$. In the rest of the section all the pullbacks and fibres are homotopy.

Combinatorics is partially about counting structures, and this is done by taking cardinality of sets. However, sometimes this does not directly give the desired result, because, as mentioned above,  symmetries between the structures have to be taken into account. These symmetries cannot be encoded inside a set, but can be encoded inside a groupoid. If a combinatorial structure is encoded in a groupoid rather than in a set, then we need a notion of cardinality for groupoids.  The \emph{homotopy cardinality} \cite[\S 3]{GKT:HLA} of a finite groupoid $X$ is defined as 
$$|X|:=\sum_{x\in \pi_0X}\frac{1}{|\aut( x)|}\in \Q.$$
Notice that if $X$ is a set, that is, $\aut( x)=1$ for every $x\in X$, then its homotopy cardinality coincides with its cardinality. For $B$ a groupoid, the homotopy objective counterpart of the vector space $\Q_{\pi_0B}$ spanned by ${\pi_0B}$ is the slice category $\infgrpd_{/B}$. A finite map of groupoids $Y\xrightarrow{p}B$ corresponds to the vector
$$|p|:=\sum_{b\in \pi_0B}\frac{|Y_b|}{|\aut(b)|}\delta_b,$$
called the \emph{homotopy cardinality} of $p$. In this sum $Y_b$ is the homotopy fibre, and $\delta_b$ is a formal symbol representing the isomorphism class of $b$. A simple computation shows that $|1\xrightarrow{\name{b}}B|=\delta_b$. 

The importance of factorizations of arrows in incidence coalgebras suggests a simplicial viewpoint. This leads to the generalization of Leroux theory to $\infty$-categories and \emph{decomposition spaces}, developed by Galvez, Kock and Tonks \cite{GKT:HLA,GKT:DSIAMI-1,GKT:DSIAMI-2}. These are a kind of simplicial spaces that express the ability to decompose. It is worth mentioning that decomposition spaces encode many more combinatorial coalgebras than merely those arising from posets or categories. In this survey, however, we will only deal with Segal spaces, a particular case of decomposition space which express the ability to compose, besides the ability to decompose.

 We denote by $\simplexcategory$ the \emph{simplex category}, whose objects are finite nonempty standard ordinals 
 $$[n]=\{0<1<\cdots<n\}$$
 and whose morphisms are order preserving maps between them. These maps are generated by the \emph{coface} maps $\partial^i\colon [n-1]\rightarrow [n]$ which skips $i$, and the \emph{codegeneracy} maps $\sigma^i\colon [n+1]\rightarrow [n]$ which repeats $i$. The obvious relations between this maps, such as $\partial^i\partial^j=\partial^{j-1}\partial^i$ for $i<j$, are called \emph{cosimplicial identities}.
 
A simplicial groupoid  is a (pseudo-)functor $X\colon \simplexcategory^{\text{op}}\longrightarrow \infgrpd$. The image of $[n]$ is denoted by $X_n$ and called the groupoid of $n$-simplices. The images of $\partial^i$ and $\sigma^i$ are denoted $d_i$ and $s_i$ and called face and degeneracy maps respectively. Explicitely, a simplicial groupoid is a sequence of groupoids $(X_n)_{n\geq 0}$ together with morphisms $d_i\colon X_n\rightarrow X_{n-1}$ and $s_i\colon X_n\rightarrow X_{n+1}$ for $0\leq i \leq n$, satisfying the \emph{simplicial identities}, induced by the cosimplicial identities, such as $d_id_j\simeq d_{j-1}d_i$ for $i<j$. The symbol $\simeq$ means that the identities are not equalities but coherent isomorphisms. This, roughly speaking, is what the (pseudo-) above represents, but we do not have to worry about this here.

A simplicial groupoid $X$ is a \emph{Segal space} \cite[\S 2.9, Lemma 2.10]{GKT:DSIAMI-1} if the following square is a pullback for all $n>0$:
\begin{equation}\label{segal}
\begin{tikzcd}[column sep=small]
X_{n+1} \arrow[r,"d_0"] \arrow[d,"d_{n+1}"'] \rdpbk& X_n \arrow[d,"d_n"]\\
	      X_n \arrow[r,"d_0"'] & X_{n-1}.
 \end{tikzcd}
\end{equation}

Segal spaces arise prominently through the fat nerve construction: the \emph{fat nerve} of a category $\mathcal{C}$ is the simplicial groupoid $X=N\mathcal{C}$ with $X_n=\fun ([n],\C)^{\simeq}$, the groupoid of functors $[n] \to \C$. 
%For example the objects and morphisms of $X_3$ can be pictured as
%\begin{center}
%\begin{tikzcd}[sep={3em,between origins}]
%a\arrow[r]&b\arrow[r]&c\arrow[r]&d
%\end{tikzcd}
%\hspace{15pt}
%\text{and}
%\hspace{15pt}
%\begin{tikzcd}[sep={3em,between origins}]
%a\arrow[r]\ar{d}{\wr}&b\arrow[r]\arrow[d]&c\arrow{d}{\wr}\arrow[r]&d\arrow{d}{\wr}\\
%a'\arrow[r]&b'\arrow[r]&c'\arrow[r]&d',
%\end{tikzcd}
%\end{center}
%where the horizontal and vertical arrows are respectively morphisms and isomorphisms in $\C$.
In this case the pullbacks above are strict, so that all the simplices are strictly determined by $X_0$ and $X_1$, respectively the objects and arrows of $\C$, and the inner face maps are given by composition of arrows in $\C$. In the general case $X_n$ is determined from $X_0$ and $X_1$ only up to equivalence, but one may still think of it as a ``category'' object whose composition is defined only up to equivalence.

Let $X$ be a simplicial groupoid. 
The spans
$$X_1\xleftarrow{\;\;d_1\;  \;} X_2 \xrightarrow{(d_2,d_0)\;}X_1\times X_1,    \;\;\;\;\;\;\;\;\;\;\;\;\;\; X_1\xleftarrow{\;\;s_0\; \;} X_0 \xrightarrow{\;\;t\;\;}1, $$

\noindent define two functors 
\begin{equation}\label{segalcom}
\begin{tabular}{rllccrll}
$\Delta \colon \infgrpd_{/X_1}$& $\longrightarrow$ &$\mathbf{Grpd}_{/X_1\times X_1}$ && & $\epsilon \colon \infgrpd_{/X_1}$& $\longrightarrow$ &$\infgrpd$\\
$S\xrightarrow{s} X_1$& $\longmapsto$ & $(d_2,d_0)_!\circ d_1^{\ast}(s)$, && &$S\xrightarrow{s} X_1$& $\longmapsto$ & $t_!\circ s_0^{\ast}(s)$ .
\end{tabular}
\end{equation} 
Recall aso that upperstar is pullback and lowershriek is postcomposition. This is the general way in which spans interpret homotopy linear algebra \cite{GKT:HLA}.

As mentioned above, Segal spaces are a particular case of decomposition spaces \cite[Proposition 3.7]{GKT:DSIAMI-1}, simplicial groupoids with the property that the functor $\Delta$ is coassociative with the functor $\epsilon$ as counit (up to homotopy). In this case $\Delta$ and $\epsilon$ endow $\infgrpd_{/X_1}$ with a coalgebra structure \cite[\S 5]{GKT:DSIAMI-1}   called the \emph{incidence coalgebra} of $X$. Note that in the special case where $X$ is the nerve of a poset, this
construction becomes the classical incidence coalgebra construction after
taking cardinality, as we shall do shortly.

A Segal space $X$ is \emph{CULF monoidal} \cite[\S 4,9]{GKT:DSIAMI-1} if it has a product $X_n\times X_n\rightarrow X_n$ for each $n$, compatible with the degeneracy and face maps, and such that for all $n$ the squares
\begin{equation}
\label{monoidalpullback}
\begin{tikzcd}
X_n\times X_n \arrow[r,"g\times g"] \arrow[d] \rdpbk & X_1\times X_1 \arrow[d]\\
	      X_n \arrow[r,"g"] & X_1.
 \end{tikzcd}
\end{equation}
where $g$ is induced by the unique endpoint-preserving map $[1]\rightarrow [n]$, are pullbacks \cite[\S  4]{GKT:DSIAMI-1}. For example the fat nerve of a monoidal extensive category is a CULF monoidal Segal space. Recall that a category $\C$ is monoidal extensive if it is monoidal $(\C,+,0)$ and the natural functors $\C_{/A}\times \C_{/B}\rightarrow \C_{/A+B}$  and $\C_{/0}\rightarrow 1$ are equivalences.
 
If $X$ is CULF monoidal then the resulting coalgebra is in fact a bialgebra \cite[\S 9]{GKT:DSIAMI-1}, with product given by
\begin{center}
\begin{tabular}{rllccrll}
$\grpdSprod\colon \infgrpd_{/X_1}\otimes \infgrpd_{/X_1}$&$\xrightarrow{\; \sim \;}$& $\infgrpd_{/X_1\times X_1}$  &$\xrightarrow{\;\; +_! \;\;}$ &$\infgrpd_{/X_1}$ \\
$(G\rightarrow X_1)\otimes (H\rightarrow X_1)$& $\longmapsto$ & $G\times H\rightarrow X_1\times X_1$&$ \longmapsto $& $ G\times H\rightarrow X_1$.
\end{tabular}
\end{center} 
Briefly, a product in $X_n$ compatible with the simplicial structure endows $X$ with a product, but in order to be compatible with the coproduct it has to satisfy the diagram \eqref{monoidalpullback} (i.e. it has to be a CULF functor).

 A Segal space $X$ is \emph{locally finite} \cite[\S 7]{GKT:DSIAMI-2} if $X_1$ is a locally finite groupoid and both $s_0\colon X_0\rightarrow X_1$ and $d_1\colon X_2\rightarrow X_1$ are finite maps. In this case one can take homotopy cardinality to get a comultiplication
\begin{center}
\begin{tabular}{rll}
$\Delta \colon \Q_{\pi_0X_1}$& $\longrightarrow$ &$ \Q_{\pi_0X_1}\otimes  \Q_{\pi_0X_1}$\\
$|S\xrightarrow{s} X_1|$& $\longmapsto$ & $|(d_2,d_0)_!\circ d_1^{\ast}(s)|$
\end{tabular}
\end{center} 
and similarly for $\epsilon$ (cf. \cite[\S 7]{GKT:DSIAMI-2}). Moreover, if $X$ is CULF monoidal then $\Q_{\pi_0X_1}$ acquires a bialgebra structure with the product $\cdot=|\grpdSprod|$. In particular, if we denote by $+$ the monoidal product in $X$, then
$\delta_{a}\cdot \delta_{b}=\delta_{a+b}$ for any $|1\xrightarrow{\name{a}}X_1|$ and $|1\xrightarrow{\name{b}}X_1|$. The following lemma gives a formula to compute the coproduct of the isomorphism class of an element $f\in X_1$.

\begin{lemma}[{\cite[\S 4]{Cebrian}}]
\label{SegalCom} 
Let $X$ be a Segal space. Then for $f$ in $X_1$ we have
$$\Delta(\delta_f)=\sum_{b\in \pi_0X_1}\sum_{a\in \pi_0X_1} \frac{|\iso(d_0a,d_1b)_{f}|}{|\aut(b)||\aut(a)|}\delta_a\otimes \delta_b.$$
\end{lemma}
\noindent Here $\iso(d_0a,d_1b)$ refers to the set of isomorphisms from $d_0a$ to $d_1b$, and the subscript $f$ means homotopy fibre.
%%----------------------------------------------------------------------------------------
%%	SECTION: species and faa di bruno
%%----------------------------------------------------------------------------------------
\section{Species and the Fa\`a di Bruno bialgebra}\label{species}

One of the starting points for objective combinatorics is the theory of species, introduced by Joyal \cite{Joyal:1981}. Through the notion of species, Joyal showed that manipulations with power series and generating functions can be carried out directly on the combinatorial structures themselves. A \emph{species} is a functor 
$$F\colon \B \longrightarrow \B$$
from the category $\B$ of finite sets and bijections to itself. To each finite set $S$ the species $F$ associates another finite set $F[S]$, whose elements are called $F$-\emph{structures} on the set $S$. Each bijection $S\rightarrow R$  gives a bijection $F[S]\rightarrow F[R]$. For example, there is a species $\Pi$ that sends a set $E$ to $\Pi(E)$, the set of all its partitions, and a bijection $E\xrightarrow{f} E'$ to the obvious bijection $\Pi(E)\rightarrow \Pi(E')$ given by $f$. Other examples include structures of graphs, trees, linear orders, etc. 

We may attach different kinds of power series to a species $F$ in order to enumerate the $F$-structures or the isomorphism classes of $F$-structures. The first ones are often called \emph{labelled} structures, while the second ones are called \emph{unlabelled} structures. The \emph{exponential generating function} associated to $F$ is 
$$F(x)=\sum_{n\geq0}|F[n]|\frac{x^n}{n!},$$
where $|F[n]|$ is the number of $F$-structures on a set of $n$ elements. This function is used for labelled enumeration. The \emph{type generating function} associated to $F$ is 
$$\tilde{F}(x)=\sum_{n\geq 0} |\tilde{F}[n]|x^n,$$
where $|\tilde{F}[n]|$ is the number of unlabelled structures of $F$. For example it is easy to see that $|\Pi[3]|=5$, while $|\tilde{\Pi}[3]|=3$, because the three partitions $\{\{1,2\},\{3\}\}, \{\{1\},\{2,3\}\}$ and $\{\{1,3\},\{2\}\}$ are isomorphic. 

Several operations of generating functions can be lifted to the level of species \cite{Joyal:1981}. For instance, given two species $F$ and $G$ we define their \emph{sum} and \emph{product} by
$$(F+G)[S]=F[S]+G[S] \;\; \text{and} \;\; (F\cdot G)[S]=\sum_{\substack{S_1+S_2=S\\S_1\cap S_2=\o}}F[S_1]\times G[S_2]$$
respectively. Both operations are compatible with sum and multiplication of generating functions, so that $(F+G)(x)=F(x)+G(x)$, $(F\cdot G)(x)=F(x)\cdot G(x)$ and similarly for the type generating functions. Nevertheless, the operation that interests us most is \emph{substitution} \cite[\S 2.2]{Joyal:1981}. Suppose that $G[\o]=\o$. Then
\begin{equation}\label{speciessubs}
(F\circ G)[S]=\sum_{\pi \in \Pi(S)} F[\pi]\times \prod_{B\in \pi} G[B],
\end{equation}
where $F[\pi]$ interprets $\pi$ as a set. Notice that this is not the composition of $F$ and $G$ as functors. Substitution of species is compatible with the exponential generating function, but not with the type generating function. To obtain a power series for unlabelled enumeration compatible with substitution a third kind of generating function is required. 
The \emph{cycle index series} of a species $F$ \cite[\S 3]{Joyal:1981} is the formal power series (in infinitely many variables $x_1,x_2,\dots$)
$$Z_F(x_1,x_2,\dots)=\sum_{n\geq 0}\frac{1}{n!}\left (\sum_{\sigma \in \mathfrak{S}_n} |\text{Fix}(F[\sigma])x_1^{\sigma_1}x_2^{\sigma_2}\cdots\right),$$
where $\mathfrak{S}_n$ denotes the group of permutations of $[n]$, $\sigma_k$ is the number of cycles of size $k$ of $\sigma$ and $\text{Fix}(F[\sigma])$ is the set of $F$-structures fixed by $F[\sigma]$.
For example, it is easy to see that 
$$\sum_{\sigma \in \mathfrak{S}_3} |\text{Fix}(\Pi[\sigma])x_1^{\sigma_1}x_2^{\sigma_2}\cdots=5x_1^3+3x_1x_2+x_3.$$
The cycle index series was first used by P\'olya \cite{polya1937}. It satisfies that $Z_F(x,0,\dots)=F(x)$ and $Z_F(x,x^2,\dots)=\tilde{F}(x)$. Furthermore, it is compatible with the following notion of substitution:
\begin{definition}[\cite{polya1937}]Given two power series, $F(x_1,x_2,\dots)$ and $G(x_1,x_2,\dots)$, their \emph{plethystic substitution} is defined as 
$$(G\oast F)(x_1,x_2,\dots)=G(F_1,F_2,\dots),$$
with $F_k=F(x_k,x_{2k},\dots)$.
\end{definition}

Back to the definition of substitution of species (\ref{speciessubs}), observe that the relevant information comes from a decomposition of $S$. This decomposition is in fact the comultiplication of the isomorphism class of the interval $[\hat{0},\hat{1}]$ of partitions of $S$ in the incidence bialgebra of the poset of partitions. Indeed, 
\begin{equation}\label{partcom}
\Delta([\hat{0},\hat{1}])=\sum_{\hat{0}\leq \pi \leq\hat{1}} [\hat{0},\pi]\otimes [\pi,\hat{1}],
\end{equation}
which  becomes the same as (\ref{speciessubs}) after expressing partitions as the disjoint union of their blocks. In fact disjoint union gives this coalgebra a structure of bialgebra, known as the \emph{Fa\`a di Bruno bialgebra}. 

In purely algebraic terms, the \emph{Fa\`a di Bruno bialgebra} $\mathcal{F}$ is the free algebra $\Q[A_1,A_2,\dots]$, where $A_n$ is the dual map $A_n\in \Q\pow{x}^{\ast}$ defined by 
$$A_n(f)=\frac{d^nf}{dx^n}.$$
Its comultiplication is defined to be dual to substitution of power series. That is
$$\Delta(A_n)(F\otimes G)=A_n(G\circ F).$$
The comultiplication of $A_n$ corresponds to the comultiplication of $[n]$ in the incidence coalgebra of partitions, and can be expressed with the \emph{Bell polynomials} $B_{n,k}(A_1,A_2,\dots)$, which count the number of partitions of a set with $n$ elements into $k$ blocks:
$$\Delta(A_n)=\sum_{k=1}^n B_{n,k}(A_1,A_2,\dots)\otimes A_k.$$
 The category of partitions is equivalent to the category of surjections, so that $\mathcal{F}$ can be expressed from surjections too \cite[\S 7.4]{Joyal:1981}, and in fact it looks simpler. In view of (ii) of Lemma \ref{lema:PS}, equation (\ref{partcom}) corresponds to 
 $$\Delta(S\twoheadrightarrow 1)=\sum_{S\twoheadrightarrow R \twoheadrightarrow 1} (S\twoheadrightarrow R) \otimes (R \twoheadrightarrow 1).$$
 The algebra structure is again given by disjoint union of sets. This sum is over isomorphism classes of factorizations $S\twoheadrightarrow R \twoheadrightarrow 1$, meaning up to isomorphism $R\xrightarrow{\sim}R'$ making the diagram commute. The precise statement that this comultiplication on surjections (or partitions) gives in fact the Fa\`a di Bruno bialgebra fits very well into the theory of Segal spaces, where all the issues with isomorphism classes take care of themselves. We denote by $\sur$ the category whose objects are finite sets and whose morphisms are surjections.
 \begin{theorem}[\cite{Joyal:1981,GKT:FdB}]
The Fa\`a di Bruno bialgebra $\mathcal{F}$ is equivalent to $\Q_{\pi_0\Sg}$, the homotopy cardinality of the incidence bialgebra of the fat nerve $N\sur\colon \simplexcategory^{\text{op}}\rightarrow \infgrpd$ of the category of surjections. 
\begin{proof}
Notice that $\Sg\simeq (N\sur)_1$. We know that $\mathcal{F}$ is generated by the functionals $A_n$. Any surjection is isomorphic to the disjoint union of surjections with singleton target. Hence,  ${\delta_n=|1\xrightarrow{\name{n\twoheadrightarrow 1}} \Sg|}$ corresponds to $A_n$. Using Lemma \ref{SegalCom} we get
$$\Delta (\delta_n)=\sum_{b:k\twoheadrightarrow 1}\sum_{a:n\twoheadrightarrow k} \frac{|\iso(k,k)_{n\twoheadrightarrow 1}|}{|\aut(b)||\aut(a)|}\delta_a\otimes \delta_k.$$
It is clear that $|\aut(b)|=k!$. Now, in this case any element of $\iso(k,k)$ gives $n\twoheadrightarrow 1$, so that $|\iso(k,k)_{n\twoheadrightarrow 1}|=n!\cdot k!$ (the $n!$ appears because we are taking homotopy cardinality).  Moreover, $\delta_a=\delta_{n_1}\cdots \delta_{n_k}$, where $n_i$ are the fibres of $a\colon n\twoheadrightarrow k$.
%, and if we denote by $\l$ the type of $a$ then, as we already know, $|\aut(a)|=\displaystyle\prod_i \l_i!\cdot i!^{\l_i}$. 
Altogether we obtain
$$\Delta (\delta_n)=\sum_{k\twoheadrightarrow 1}\sum_{n\twoheadrightarrow k} \frac{n!}{|\aut(n\twoheadrightarrow k)|}\prod_{i=1}^k \delta_{n_i}\otimes \delta_k,$$
which is easily checked to correspond to the comultiplication of $A_n$.
\end{proof}
\end{theorem}

We would like to give an  interpretation of plethystic substitution from partitions and surjections, as we have just done for substitution of ordinary power series. The monomials of ordinary power series are indexed by natural numbers, which coincide with isomorphism classes of sets. However, the monomials of power series in infinitely many variables are indexed by isomorphism classes of partitions $\l=(\l_1,\l_2,\dots)$. This is why Nava and Rota \cite{Nava-Rota} developed the notion of partitional as a generalization of species, in order to give an interpretation of plethystic substitution analogous to the species interpretation of ordinary substitution. A sequence like $\l$ could also represent the isomorphism class of a permutation, and in fact Bergeron \cite{Bergeron} gave a similar interpretation but in terms of permutationals, rather than partitionals. 

 %%----------------------------------------------------------------------------------------
%%	SECTION: partitionals and plethystic substitution
%%----------------------------------------------------------------------------------------
\section{Partitionals and the plethystic bialgebra}\label{partitionals}

A partitional \cite{Nava-Rota} is a functor $M\colon \bbP\rightarrow \B$ from the category of partitions $\bbP$ to the category of sets ans bijections $\B$. The image $M[E,\pi]$ of $(E,\pi)$ under $M$ is the set of $M$-\emph{structures}. By functoriality, the cardinality $|M[E,\pi]|$ depends only on the isomorphism class $\l$ of the partition $(E,\pi)$, and will be denoted by $M[\l]$. Therefore we can define the \emph{generating function} of $M$ as
\begin{equation}\label{MoR}
M(x_1,x_2,\dots)=\sum_{\l}\frac{M[\l]}{\autiv{\l}}x_1^{\l_1}x_2^{\l_2}\cdots.
\end{equation}

As in the case of species, several operations of generating functions can be lifted to the level of partitionals. The sum and the product are defined in a similar way as for species, and we will not do it here. The substitution, however, is more complex and involves the notion of transversal. Let $M$ and $R$ be two partitionals, then their \emph{substitution} \cite[\S 6]{Nava-Rota} is defined as 
$$(M\circ R)[E,\sigma]:=\sum_{\substack{(\pi,\tau)\\ \text{transversal of } \sigma}} M[\tau,(\sigma\vee\tau)|\tau]\times \prod_{B\in \sigma \vee \tau}R[\pi_B,\sigma_B|\pi_B].$$
Substitution of partitionals is compatible with plethystic substitution of generating functions. That is, $(M\circ R)(x_1,x_2,\dots)=M(x_1,x_2,\dots)\oast R(x_1,x_2,\dots)$ \cite[\S 6]{Nava-Rota}. Notice that, as before, this definition is also based on a decomposition of $(E,\sigma)$, which under isomorphism classes gives rise to a comultiplication in a bialgebra, the \emph{plethystic bialgebra} \cite[\S 3]{Cebrian}. 
For each $\lambda$ define the functional $A_{\l}\in \Q\pow{\x}^{\ast}$ by  $A_{\l}(F)=f_{\l}$. The \emph{plethystic bialgebra} $\mathcal{P}$ is the free polynomial algebra $\Q[\{A_{\l}\}_{\l}]$ along with the comultiplication dual to plethystic substitution. That is, for each $\lambda$ and $F,G \in \Q\pow{\x}$,
$$ \Delta(A_{\l})(F\otimes G)=A_{\l}(G\oast F).$$
The counit is given by $\epsilon (A_{\l})=\langle A_{\l}, x_1\rangle$. 

It is difficult to express $\mathcal{P}$ as an incidence coalgebra from partitions, but using the machinery of Section \ref{preliminaries} we can find a Segal groupoid, arising also from surjections, whose incidence bialgebra is isomorphic to the plethystic bialgebra. First of all, let us apply the results of Section \ref{partitions} to express definition (\ref{MoR}) in the context of surjections:
\begin{equation}\label{comTS}
\Delta\left(\begin{tikzcd}[sep={2em,between origins}]
& E \arrow[ld,twoheadrightarrow] \arrow[rd,twoheadrightarrow]  \arrow[d,phantom,"\sigma"]& \\
             B \arrow[rr,twoheadrightarrow]& \phantom{a} &1
\end{tikzcd}\right)
:=\hspace{-40pt}\sum_{
\begin{tikzcd}[sep={2em,between origins}]
& & E \arrow[ld,twoheadrightarrow] \arrow[rd,twoheadrightarrow]\dpbk& & \\                             
			&S\arrow[ld,twoheadrightarrow] \arrow[rd,twoheadrightarrow]  \arrow[d,phantom,"\pmb{\mu}"]& & X\arrow[ld,twoheadrightarrow] \arrow[rd,twoheadrightarrow] \arrow[d,phantom,"\l"]&\\
	      B\arrow[rr,twoheadrightarrow]& \phantom{a}&I \arrow[rr,twoheadrightarrow]& \phantom{a}&1
\end{tikzcd}}
\hspace{-40pt}
 \begin{tikzcd}[sep={2em,between origins}]
& S \arrow[ld,twoheadrightarrow] \arrow[rd,twoheadrightarrow]  \arrow[d,phantom,"\pmb{\mu}"]& \\
             B \arrow[rr,twoheadrightarrow]&\phantom{a}  &I
\end{tikzcd}
\otimes
\begin{tikzcd}[sep={2em,between origins}]
& X \arrow[ld,twoheadrightarrow] \arrow[rd,twoheadrightarrow]  \arrow[d,phantom,"\l"]& \\
             I \arrow[rr,twoheadrightarrow]& \phantom{a} &1.
\end{tikzcd}
 \end{equation}
 
\noindent Again, this sum is over isomorphism classes of transversals of surjections, meaning up to isomorphisms $S\xrightarrow{\sim} S'$ and $X\xrightarrow{\sim} X'$ making the diagram commute. To give a precise statement about the relation between surjections and the plethystic bialgebra we introduce the simplicial groupoid $T\sur$ \cite[\S 2]{Cebrian}, whose $n$-simplices are pyramids like the one pictured below for $n=3$
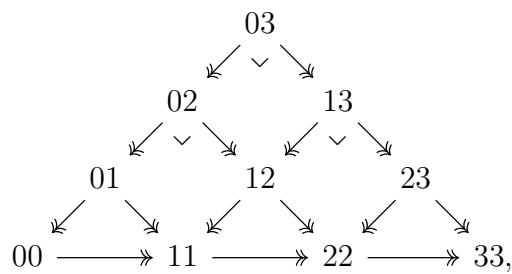
\begin{center}
\begin{tikzcd}[sep={2.5em,between origins}]
& & &{03} \arrow[ld,twoheadrightarrow] \arrow[rd,twoheadrightarrow]\dpbk& &&\\
		& & {02} \arrow[ld,twoheadrightarrow] \arrow[rd,twoheadrightarrow]\dpbk& & {13} \arrow[ld,twoheadrightarrow] \arrow[rd,twoheadrightarrow]\dpbk& &\\                             
			&{01}\arrow[ld,twoheadrightarrow] \arrow[rd,twoheadrightarrow]& &{12}\arrow[ld,twoheadrightarrow] \arrow[rd,twoheadrightarrow] & & 23 \arrow[ld,twoheadrightarrow] \arrow[rd,twoheadrightarrow]&\\
	     {00} \arrow[rr,twoheadrightarrow]& &{11} \arrow[rr,twoheadrightarrow]& &{22} \arrow[rr,twoheadrightarrow]& &{33}, 
\end{tikzcd}
\end{center}
where all the entries are finite sets, the arrows are surjections, and the squares are pullbacks of sets. Morphisms in $T\sur_n$ are levelwise bijections making the diagram commute. The face map $d_i$ removes all the object containing an $i$ index, while the degeneracy maps $s_i$ repeats the objects containing an $i$ index. The fact that the squares are pullbacks makes it a Segal space. For instance, 
\begin{center}
\begin{tikzcd}[sep={2em,between origins}]                 
			&{01}\arrow[ld,twoheadrightarrow] \arrow[rd,twoheadrightarrow]& \arrow[rd,twoheadrightarrow,phantom,"\circ"]& & 12 \arrow[ld,twoheadrightarrow] \arrow[rd,twoheadrightarrow]&\\
	     {00} \arrow[rr,twoheadrightarrow]& &{11}  &{11} \arrow[rr,twoheadrightarrow]& &{22}, 
\end{tikzcd}
$=d_1 \left(
\begin{tikzcd}[sep={2em,between origins}]
		& & {02} \arrow[ld,twoheadrightarrow] \arrow[rd,twoheadrightarrow]\dpbk& &\\                             
			&{01}\arrow[ld,twoheadrightarrow] \arrow[rd,twoheadrightarrow]& &{12}\arrow[ld,twoheadrightarrow] \arrow[rd,twoheadrightarrow] &\\
	     {00} \arrow[rr,twoheadrightarrow]& &{11} \arrow[rr,twoheadrightarrow]& &{22}
\end{tikzcd}\right)
=
\begin{tikzcd}[sep={2em,between origins}]                           
			&{02}\arrow[ld,twoheadrightarrow] \arrow[rd,twoheadrightarrow]& \\
	     {00} \arrow[rr,twoheadrightarrow]& &{22},
\end{tikzcd}$
\end{center}
which is defined up to isomorphism. Also, disjoint union of diagrams makes it a CULF monoidal Segal space \cite[\S 2]{Cebrian}.
\begin{theorem}[{\cite[\S 4]{Cebrian}}] The plethystic bialgebra $\bbP$ is isomorphic to $\Q_{T\sur_1}$, the homotopy cardinality of the incidence bialgebra of $T\sur$. 
\end{theorem}
Proving this result would be beyond the scope of this paper. The idea is that $A_{\sigma}$ corresponds to
$$\delta_{\sigma}=|1\xrightarrow{\name{\sigma}} T\sur_1|.$$
The homotopy fibre 
\begin{center}
\begin{tikzcd}
(T\sur_1)_\sigma \arrow[r]\arrow[d]\rdpbk& T\sur_2\arrow[d,"d_1"]\\
1\arrow[r,"\name{\sigma}"]& T\sur_1
\end{tikzcd}
\end{center}
is precisely the groupoid of factorizations appearing in the subindex of the comultiplication (\ref{comTS}) and the summands come respectively from $d_2$ and $d_0$ of the factorizations.

%%----------------------------------------------------------------------------------------
%%	REFERENCES
%%----------------------------------------------------------------------------------------

%%----------------------------------------------------------------------------------------
%%	END
%%----------------------------------------------------------------------------------------


\begin{thebibliography}{13}


%\bibitem{bauer}
%{\sc Tilman Bauer}.
%\newblock {\em Formal plethories}.
%\newblock Adv. Math. {\bf 254} (2014), 497--569.
%\newblock \arxiv{1107.5745}.

\bibitem{Bergeron}
{\sc Fran\c{c}ois Bergeron}.
\newblock {\em Une combinatoire du pl\'ethysme}.
\newblock J. Combin. Theory Ser. A {\bf 46} (1987), 291--305.

\bibitem{Cebrian}
\newblock{\sc Alex Cebrian}.
\newblock{\em A simplicial groupoid for plethysm}.
\newblock \arxiv{1804.09462} (2018).

\bibitem{DK}
\newblock{\sc Tobias Dyckerhoff {\rm and} Mikhail Kapranov}.
\newblock{\em Higher Segal Spaces}.
\newblock{Lecture Notes in Mathematics {\bf 2244} (2019)}. 
\newblock \arxiv{1212.3563}.

%\bibitem{Bergeron:2}
%{\sc Fran\c{c}ois Bergeron}.
%\newblock {\em A combinatorial outlook on symmetric functions}.
%\newblock J. Combin. Theory Ser. A {\bf 50} (1989), 226--234.

%\bibitem{Borger-Wieland}
%{\sc James Borger {\rm and }Ben Wieland}.
%\newblock {\em Plethystic algebra}.
%\newblock Adv. Math. {\bf 194} (2005), 246--283.
%\newblock \arxiv{math/0407227}.

%\bibitem{BFK:0406117}
%{\sc Christian Brouder, Alessandra Frabetti, {\rm and }Christian
%  Krattenthaler}.
%\newblock {\em Non-commutative {H}opf algebra of formal diffeomorphisms}.
%\newblock Adv. Math. {\bf 200} (2006), 479--524.
%\newblock \arxiv{math/0406117}.

\bibitem{GKT:FdB}
{\sc Imma G{\'a}lvez-Carrillo, Joachim Kock, {\rm and }Andrew Tonks}.
\newblock {\em Groupoids and {F}a{\`a} di {B}runo formulae for {G}reen
  functions in bialgebras of trees}.
\newblock Adv. Math. {\bf 254} (2014), 79--117.
\newblock \arxiv{1207.6404}.

\bibitem{GKT:HLA}
{\sc Imma G{\'a}lvez-Carrillo, Joachim Kock, {\rm and }Andrew Tonks}.
\newblock {\em Homotopy linear algebra}.
\newblock Proc. Royal Soc. Edinburgh A. {\bf 148} (2018), 293--325.
\newblock \arxiv{1602.05082}.

\bibitem{GKT:DSIAMI-1}
{\sc Imma G{\'a}lvez-Carrillo, Joachim Kock, {\rm and }Andrew Tonks}.
\newblock {\em Decomposition spaces, incidence algebras and {M}{\"o}bius
  inversion I: basic theory}.
\newblock Adv. Math. {\bf 331} (2018), 952--105.
\newblock \arxiv{1512.07573}.

\bibitem{GKT:DSIAMI-2}
{\sc Imma G{\'a}lvez-Carrillo, Joachim Kock, {\rm and }Andrew Tonks}.
\newblock {\em Decomposition spaces, incidence algebras and {M}{\"o}bius
  inversion II: completeness, length filtration, and finiteness}.
  \newblock Adv. Math. {\bf 333} (2018), 1242--1292.
\newblock \arxiv{1512.07577}.



\bibitem{Joyal:1981}
{\sc Andr{\'e} Joyal}.
\newblock {\em Une th\'eorie combinatoire des s\'eries formelles}.
\newblock Adv. Math. {\bf 42} (1981), 1--82.

\bibitem{LawvereMenni}
\newblock{\sc  F. William Lawvere and Mat\'ias Menni}.
\newblock{\em The Hopf algebra of M{\"o}bius intervals.} 
\newblock Theory Appl. Categ. {\bf 24} (2010), 221--265.

%\bibitem{Littlewood305}
%{\sc Dudley E. Littlewood}.
%\newblock {\em Invariant theory, tensors and group characters}.
%\newblock Philosophical Transactions of the Royal Society of London A:
%  Mathematical, Physical and Engineering Sciences {\bf 239} (1944), 305--365.

%\bibitem{Lurie:HA}
%{\sc Jacob Lurie}.
%\newblock {\em Higher Algebra}.
%\newblock Available from \url{http://www.math.harvard.edu/~lurie/}, 2013.
\bibitem{Leroux}
\newblock{\sc Pierre Leroux}.
\newblock{\em Les cat\'egories de M{\"o}bius}.
\newblock Cahiers Topol. G\'eom. Diff. {\bf 16} (1975), 280--282.

\bibitem{Nava-Rota}
{\sc Oscar Nava {\rm and }Gian-Carlo Rota}.
\newblock {\em Plethysm, categories, and combinatorics}.
\newblock Adv. Math. {\bf 58} (1985), 61--88.

\bibitem{polya1937}
{\sc George P{\'o}lya}.
\newblock {\em Kombinatorische Anzahlbestimmungen f{\"u}r Gruppen, Graphen und
  chemische Verbindungen}.
\newblock Acta Math. {\bf 68} (1937), 145--254.

\bibitem{Rota} 
{\sc Gian-Carlo Rota}.
\newblock {\em On the foundations of combinatorial theory I. Theory of M{\"o}bius functions}.
\newblock Z. Wahrscheinlichkeitstheorie und Verw. Gebiete. {\bf 2} (1964), 340--368.

\end{thebibliography}
 \end{document}